\documentclass{amsproc}
\usepackage{amssymb}

\newtheorem{theorem}{Theorem}[section]

\newtheorem{conjecture}[theorem]{Conjecture}

\newtheorem*{DC}{Dyson's Conjecture}
\newtheorem*{qDC}{Andrews' $q$-Dyson Conjecture}

\theoremstyle{definition}

\theoremstyle{remark}

\numberwithin{equation}{section}

\newcommand{\xv}{\mathbf{x}}

\newcommand{\av}{\mathbf{a}}

\newcommand{\s}{\sigma_n(\av)}

\begin{document}

\title{Disturbing the $q$-Dyson Conjecture}

\author{Andrew V. Sills}
\address{Department of Mathematical Sciences\\ 
Georgia Southern Univeristy\\ Statesboro, GA 30460}
\email{ASills@GeorgiaSouthern.edu}

%    General info
%\date{\today}
\date{December 7, 2007}


\begin{abstract}
I discuss the computational methods behind the formulation of some conjectures
related to variants on Andrews' $q$-Dyson conjecture.
\end{abstract}

\maketitle

\section{Introduction}
  In 1962~\cite{D2}, Freeman Dyson made the following conjecture:
\begin{DC}  For positive integers $n$  and $a_1, a_2, \dots, a_n$,
  the constant term in the expansion of the Laurent polynomial
  \begin{equation} \label{DysonProd}
  \prod_{1\leqq i < j \leqq n} \left ( 1 - \frac{x_i}{x_j} \right)^{a_j} \left ( 1 - \frac{x_j}{x_i} \right)^{a_i} 
  \end{equation}
is the multinomial co\"efficient
\[ \frac{ (a_1 + a_2 + \cdots + a_n)! }{ a_1! a_2! \cdots a_n! }. \]
\end{DC}
  Dyson's conjecture was settled independently by Gunson~\cite{Gu} and Wilson~\cite{W}.
In 1970, Good~\cite{Go} supplied a particularly compact and elegant proof.

  In 1975, George Andrews conjectured a $q$-analog of Dyson's conjecture~\cite{A}:
\begin{qDC}  For nonnegative integers $n$  and $a_1, a_2, \dots, a_n$,
  the constant term in the expansion of the Laurent polynomial
  \[ \prod_{1\leqq i<j \leqq n}  \left( 1- \frac{x_i}{x_j} q \right)  \left( 1- \frac{x_i}{x_j} q^2 \right)
  \cdots  \left( 1- \frac{x_i}{x_j} q^{a_j} \right)   \left( 1- \frac{x_j}{x_i} \right) \left( 1- \frac{x_j}{x_i} q\right) 
  \cdots \left( 1- \frac{x_j}{x_i} q^{a_i -1 } \right)\]
is the $q$-multinomial co\"efficient
\[ \frac{ [a_1 + a_2 + \cdots + a_n]_q! }{ [a_1]_q! [a_2]_q! \cdots [a_n]_q! }, \]
\end{qDC}
where 
\[ [a]_q := \frac{1-q^{a}}{1-q} = 1+q+q^2+\cdots+q^{a-1} \]
is the usual $q$-analog of the nonnegative integer $a$ and
\[ [a]_q ! = \prod_{j=1}^a [j]_q \] is the $q$-factorial.
  Clearly, the $q=1$ case of the $q$-Dyson conjecture is the original Dyson conjecture.
The $q$-Dyson conjecture remained unsettled for a decade until it was proved
by Zeilberger and Bressoud~\cite{ZB}.  Two additional decades passed before a shorter
proof was found by Gessel and Xin~\cite{GX}.

  In~\cite{SZ}, Zeilberger and I set out to ``disturb" the Dyson conjecture\footnote{In the interest
of full disclosure, it is my esteemed coauthor for~\cite{SZ} who deserves full credit for the double pun in
our title based on~\cite{D2} and~\cite{Go}.} 
by programming
the computer to conjecture, and then provide proofs modeled after Good's proof~\cite{Go},
for closed form expressions of coefficients of terms in the expansion of~\eqref{DysonProd} \emph{other}
than the constant term.   Using our Maple package \texttt{GoodDyson}, available for 
free download from our home pages~\cite{SZ2}, the computer can 
(up to the limits imposed by time and
memory) conjecture and prove a closed form expression for
the coefficient of $x_1^{b_2} x_2^{b_2} \cdots x_n^{b_n}$ in the expansion of~\eqref{DysonProd}
for any \emph{fixed} $n$ and any \emph{fixed} $b_1, b_2, \dots, b_n$.  

  At this point, we should introduce some more notation.
For $n$ a positive integer, we define the following symbols:
  \begin{align}
     &{\av} := \langle a_1, a_2, \dots, a_n \rangle, 
         \tag{$n$-vector of symbolic nonnegative integers} \\
     &{\xv} := \langle x_1, x_2, \dots, x_n \rangle, 
          \tag{$n$-vector of indeterminants}\\
     &\s := a_1 + a_2 + \cdots + a_n, 
     \tag{first elementary symmetric polynomial in $n$ indeterminants}\\
     &(A;q)_n := \prod_{i=0}^{n-1} (1-Aq^i), \tag{rising $q$-factorial}\\
     &F_n(\xv ; \av) := \prod_{1\leqq i<j \leqq n} 
   \left(1-\frac{x_i}{x_j} \right)^{a_j}  \left(1-\frac{x_j}{x_i}\right)^{a_i},
       \tag{Dyson product}\\
     &\mathcal{F}_n( \xv ; \av ; q):= 
        \prod_{1\leqq i<j \leqq n} 
        \left( \frac{x_i q}{x_j};q \right)_{a_j}  
        \left(\frac{x_j}{x_i};q \right)_{a_i}, \label{qDysonProd}
  \tag{$q$-Dyson product}
   \end{align}
and let  $[Y]Z $ denote the coefficient of $Y$ in the expression $Z$, thus
the Dyson conjecture is
\[ [x_1^0 x_2^0 \cdots x_n^0 ] F_n(\xv; \av) = \frac{ \s !}{ a_1! a_2! \cdots a_n!}, \]
while the $q$-Dyson conjecture is
  \[ [x_1^0 x_2^0 \cdots x_n^0 ] \mathcal{F}_n(\xv; \av; q) = \frac{ [\s]_q !}{ [a_1]_q! [a_2]_q! \cdots 
  [a_n]_q!}. \]
   Using the output from many applications of the \texttt{GoodDyson} program for
various values of $n$ and $b_1, \dots, b_n$, I was able to conjecture and prove
the following ``disturbed" versions of the Dyson conjecture~\cite{Si1}:
 \begin{theorem}\label{1m1} Let $r$ and $s$ be fixed integers with 
$1\leqq r\neq s \leqq n$ and $n\geqq 2$.  Then
 \[ \left[ \frac{x_r}{x_s} \right] F_n(\xv ; \av) 
 = -\left(\frac{a_s}{1+\s -a_s}\right)
 \frac{\s!}{a_1! a_2!\cdots a_n!}. \]
 \end{theorem}
 \begin{theorem}\label{2m1m1}
  Let $r$, $s$, and $t$ be distinct fixed integers with 
$1\leqq r, s, t \leqq n$ and $n\geqq 3$.  Then
 \begin{multline*}
  \left[ \frac{x_r^2}{x_s x_t} \right] F_n(\xv ; \av) \\ = 
 \left(\frac{a_s a_t \Big( (1+\s) + (1+\s-a_s-a_t) \Big) }
 {(1+\s-a_s-a_t)(1+\s-a_s)(1+\s-a_t)}\right)
 \frac{\s!}{a_1! a_2!\cdots a_n!}.
 \end{multline*}
 \end{theorem}
\begin{theorem}\label{11m1m1}
  Let $r$, $s$, $t$, and $u$ be distinct fixed integers with 
$1\leqq r, s, t, u \leqq n$ and $n\geqq 4$.
  Then
 \begin{multline*}
  \left[ \frac{x_r x_s}{x_t x_u} \right] F_n(\xv ; \av) \\ =
 \left(\frac{a_t a_u \Big( (1+\s) + (1+\s-a_t-a_u) \Big) }
 {(1+\s-a_t-a_u)(1+\s-a_t)(1+\s-a_u)}\right)
 \frac{\s!}{a_1! a_2!\cdots a_n!}. 
 \end{multline*}
 \end{theorem}

\section{$q$-analogs of Theorems~\ref{1m1}--\ref{11m1m1}}\label{q-analogs}
Given that the Dyson conjecture has such a natural $q$-analog, it seemed
reasonable to look for comparable $q$-analogs of 
Theorems~\ref{1m1}--\ref{11m1m1}.
\subsection{Statements of the conjectures}
\begin{conjecture}[$q$-analog of Theorem~\ref{1m1}]\label{q1m1}
 Let $r$ and $s$ be fixed integers with $1\leqq r\neq s \leqq n$
and $n\geqq 2$.  Then
 \[ [x_r/x_s] \mathcal{F}_n(\xv ; \av ; q) = 
 -q^{L(r,s)}\left( \frac{[a_s]_q}{[1+\s-a_s]_q}\right)
 \frac{[\s]_q!}{[a_1]_q! [a_2]_q!\cdots [a_n]_q! }, \]
 where 
   \[ L(r,s) = \left\{ \begin{array}{ll}
                      1+\s-\sum_{k=r}^s a_k, &\mbox{if $r<s$}\\
                       \sum_{k=s+1}^{r-1} a_k, &\mbox{if $r>s$.}
                    \end{array} \right.\]
 \end{conjecture}
\begin{conjecture}[$q$-analog of Theorem~\ref{2m1m1}]\label{q2m1m1}
  Let $r$, $s$, and $t$ be distinct fixed integers with
$1\leqq r, s, t \leqq n$ and $n\geqq 3$.  Without loss of
generality we may assume that $s<t$.  Then
 \begin{multline*}
\left[ \frac{x_r^2}{x_s x_t} \right] \mathcal{F}_n(\xv ; \av ; q)\\ =
 q^{L(r,s,t)} 
 \left(\frac{[a_s]_q [a_t]_q \Big( [1+\s]_q + 
q^{M(r,s,t)}[1+\s-a_s-a_t]_q 
\Big) }
{ [1+\s-a_s-a_t]_q [1+\s-a_s]_q [1+\s-a_t]_q }\right) \\ \times
 \frac{[\s]_q!}{[a_1]_q! [a_2]_q!\cdots [a_n]_q! },
\end{multline*} 
where
   \[ L(r,s,t) = \left\{ \begin{array}{ll}
          2+2\s-2\sum_{k=r}^t a_k+\sum_{k=s+1}^{t-1}a_k, &\mbox{if $r<s<t$,}\\
          1+\s-\sum_{k=s}^t a_k + 2\sum_{k=s+1}^{r-1}a_k,&\mbox{if $s<r<t$,}\\   
          2\sum_{k=t+1}^{r-1} a_k+\sum_{k=s+1}^{t-1} a_k,&\mbox{if $s<t<r$,}
                    \end{array} \right.\] and
   \[ M(r,s,t) = \left\{ \begin{array}{ll}
            a_t, &\mbox{if $r<s<t$ or $s<t<r$,}\\
            a_s, &\mbox{if $s<r<t$.}
                      \end{array}\right.\]
 \end{conjecture}
 \begin{conjecture}[$q$-analog of Theorem~\ref{11m1m1}]\label{q11m1m1}
  Let $r$, $s$, $t$ and $u$ be distinct fixed integers with
$1\leqq r, s, t, u \leqq n$ and $n\geqq 4$.  Without loss of
generality we may assume that $r<s$ and $t<u$.  Then
 \begin{multline*}
\left[ \frac{x_r x_s}{x_t x_u} \right] \mathcal{F}_n(\xv ; \av ; q) \\ =
 q^{L(r,s,t,u)}
 \left(\frac{[a_t]_q [a_u]_q \Big( [1+\s]_q) +
q^{M(r,s,t,u)}[{1+\s-a_t-a_u}]_q
\Big) }
 {[1+\s-a_t-a_u]_q [1+\s-a_t]_q [1+\s-a_u]_q }\right) \\ \times
 \frac{[\s]_q!}{[a_1]_q! [a_2]_q!\cdots [a_n]_q! },
\end{multline*}
where
   \begin{multline*}
    L(r,s,t,u) \\ = \left\{ \begin{array}{ll}
     2+2\s-2\sum_{k=r}^u a_k+\sum_{k=r}^{s-1}a_k+\sum_{k=t+1}^{u-1} a_k, 
                        &\mbox{if $r<s<t<u$,}\\
          1+\s-\sum_{k=r}^{u} a_k+\sum_{k=t+1}^{s-1}a_k , 
                      &\mbox{if $r<t<s<u$,}\\
          1+\s-\sum_{k=r}^{s-1} a_k + 2\sum_{k=t+1}^{r-1} a_k  + 
               \sum_{k=t+1}^{u-1} a_k \\ \qquad +2\sum_{k=u+1}^{s-1} a_k ,
                       &\mbox{if $r<t<u<s$,}\\
          1+\s-\sum_{k=t}^{u} a_k+ \sum_{k=r}^{s-1} a_k 
           +2\sum_{k=t+1}^{r-1} a_k,
                    &\mbox{if $t<r<s<u$,}\\
          \sum_{k=t+1}^{r-1} a_k + \sum_{k=u+1}^{s-1}a_k, &\mbox{if $t<r<u<s$,}\\
          \sum_{k=r}^{s-1}a_k+\sum_{k=t+1}^{u-1} a_k+ 2\sum_{k=u+1}^{r-1}a_k,
                        &\mbox{if $t<u<r<s$,}
                    \end{array} \right. \end{multline*} and
   \[ M(r,s,t,u) = \left\{ \begin{array}{ll}
            a_u,         &\mbox{if $r<s<t<u$ or $r<t<u<s$ or $t<u<r<s$,}\\
            1+\s  &\mbox{if $r<t<s<u$ or $t<r<u<s$,}\\
            a_t,  &\mbox{if $t<r<s<u$.}\\
                      \end{array}\right.\]
 \end{conjecture}

\subsection{How the conjectures were formed}
Conjecture~\ref{q1m1} was found first since it is the simplest.  It is
straightforward to program a Maple procedure which extracts the
coefficient of $x_r/x_s$ of
$\mathcal{F}_n (\xv; \av; q)$
for specific 
values of $n$, $r$, $s$, $a_1, a_2, \dots, a_n$, and to divide
out the multinomial coefficient from the resulting $q$-expression.
Furthermore, the \texttt{qfactor} procedure in Frank Garvan's 
\texttt{qseries.m}
Maple package~\cite{Ga} was helpful for putting the result in a tractable form.
For $n=3$ and $4$ and various small values of $a_1, \dots, a_n$,
it became clear that to move from Theorem~\ref{1m1} to its $q$-analog,
all that was necessary was to replace each factor $z$ by $[z]_q$,
and multiply the resulting expression by $q^L$, where $L$ was an
(as yet unknown) function of the $a_i$'s that depended on $r$ and $s$.  
Upon examining the data, 
I was led to the working hypothesis
that $L$ was \emph{piecewise linear} in the $a_i$'s with different pieces
arising from some condition on $r$ and $s$.

At this point, I began to create the ``\texttt{qDysonConj}" package~\cite{Si2}.
I programmed the ``\texttt{Conj1m1}" procedure in Maple, 
which takes
as input the ordered pair $(r,s)$ and $n$, and finds the
linear function 
  \[ L(r,s) = \lambda_0 + \sum_{i=1}^n \lambda_i a_i \]
which fits the internally generated data.  

  (The name ``Conj1m1" is meant
to suggest that we wish to conjecture the missing exponent $L$ for the
coefficient of $x_1^{b_1} x_2^{b_2} \cdots x_n^{b_n}$ in the $q$-Dyson 
product (for a specific $n$) where one of the $b_i$ is 1, one of the
$b_i$ is $-1$ and the rest are zero.)  

  The idea behind the \texttt{Conj1m1} is quite simple.  Based on the
assumption
  \[ \left[ \frac{x_r}{ x_s} \right] \mathcal F_n(\xv; \av; q) =  
  -q^{ \lambda_0 + a_1 \lambda_1 + a_2 \lambda_2 + \cdots + a_n \lambda_n}
  \left( \frac{[a_s]_q}{ [1+\s-a_s]_q }\right)\frac{[\s]_q!}
  {[a_1]_q! \cdots [a_n]_q!} ,
  \]
the \texttt{Conj1m1} procedure, for a given $r$, $s$, and $n$  effectively computes
  \[  \log_q \left( \frac{ [x_r/x_s] \mathcal F_n (\xv; \av; q)}
  {-\left( \frac{[a_s]_q}{ [1+\s-a_s]_q }\right)\frac{[\s]_q!}
  {[a_1]_q! \cdots [a_n]_q!}} \right)
    \] for
$n+1$ linearly independent values of the vector $\av$
and solves the resulting linear system for $\lambda_0, \lambda_1, \dots, \lambda_{n}$.
 
   Let us recreate a Maple session to guess $L(r,s)$ using 
 the case $n=6$.

\begin{verbatim}
>  read "qDysonConj";
            Generalized qDyson conjecture package
                         by A.V. Sills
             Enter 'ez()' for a list of procedures
> C:=combinat[permute](6,2);
C := [[1, 2], [1, 3], [1, 4], [1, 5], [1, 6], [2, 1], [2, 3], [2, 4], 
[2, 5], [2, 6], [3, 1], [3, 2], [3, 4], [3, 5], [3, 6], [4, 1], [4, 2], 
[4, 3], [4, 5], [4, 6], [5, 1], [5, 2], [5, 3], [5, 4], [5, 6], [6, 1], 
[6, 2], [6, 3], [6, 4], [6, 5]]
\end{verbatim}
\begin{quote}In order to have Maple run the 
\texttt{Conj1m1} procedure on all ordered pairs $(r,s)$ with
$1\leqq r\neq s \leqq 6$, we use the built-in
\texttt{combinat[permute]} procedure, and have
Maple loop through all 30 permutations of length $2$
on the set $\{ 1, 2, 3, 4, 5, 6 \}$.
\end{quote}

\begin{verbatim}
> for k from 1 to nops(C) do Conj1m1( op(C[k]), 6) od;
                  [1, 2], 1 + a3 + a4 + a5 + a6
                    [1, 3], 1 + a4 + a5 + a6
                      [1, 4], 1 + a5 + a6
                         [1, 5], 1 + a6
                           [1, 6], 1
                           [2, 1], 0
                 [2, 3], 1 + a1 + a4 + a5 + a6
                    [2, 4], 1 + a1 + a5 + a6
                      [2, 5], 1 + a1 + a6
                         [2, 6], 1 + a1
                           [3, 1], a2
                           [3, 2], 0
                 [3, 4], 1 + a1 + a2 + a5 + a6
                    [3, 5], 1 + a1 + a2 + a6
                      [3, 6], 1 + a1 + a2
                        [4, 1], a2 + a3
                           [4, 2], a3
                           [4, 3], 0
                 [4, 5], 1 + a1 + a2 + a3 + a6
                    [4, 6], 1 + a1 + a2 + a3
                      [5, 1], a2 + a3 + a4
                        [5, 2], a3 + a4
                           [5, 3], a4
                           [5, 4], 0
                 [5, 6], 1 + a1 + a2 + a3 + a4
                   [6, 1], a2 + a3 + a4 + a5
                      [6, 2], a3 + a4 + a5
                        [6, 3], a4 + a5
                           [6, 4], a5
                           [6, 5], 0
 
\end{verbatim}

The above output shows the conjectured form of $q^{L(r,s)}$ for each
of the thirty possible values of $(r,s)$ in the case $n=6$.
Notice that when $r<s$, 
  \[ L(r,s) =  1 + \sum_{i\in\{1,2,3,4,5,6 \} \setminus \{ r, r+1, \dots, s \} } a_i,\]
while if $r>s$, 
  \[ L(r,s) = \sum_{i\in \{1,2,3,4,5,6 \} \setminus \{ s, s+1, \dots, r \}} a_i .\]

The data above, combined with the analogous data for many different values
of $n$ led me to conjecture $L(r,s)$ as given in 
Conjecture~\ref{q1m1}.

 Once I had Conjecture~\ref{q1m1}, it seemed reasonable to guess that
the $q$-analog of Theorem~\ref{2m1m1} would have an analogous form,
noting that this time the expression broke down neatly into a sum
of two terms.  I guessed that each of the two terms included a 
factor of the form $q^L$, where again, $L$ is a piecewise linear function of the
$a_i$'s that depended on $r$, $s$, and $t$; piecewise according to
the ordering of $r$, $s$, and $t$ from smallest to largest.   
This time I programmed the \texttt{Conj2m1m1} procedure which works
similarly to the \texttt{Conj1m1} procedure except that now two 
piecewise linear functions must be found simultaneously.  By
extracting the coefficient of $x^2_r/x_s x_t$ from the expanded
$q$-Dyson product for a given $r$, $s$, $t$, and $n$, and dividing
through by 
  \[ 
 \left(\frac{[a_s]_q [a_t]_q }
 { (1-q) [1+\s-a_s-a_t]_q  [1+\s-a_s]_q [1+\s-a_t]_ q }\right) 
\frac{[\s]_q!}{[a_1]_q! \cdots [a_n]_q!},\] 
what remains is a four-term polynomial in $q$ which we presume to be
of the form
   \[ q^L(1-q^{1+\s}) + q^M(1-q^{1+\s-a_s-a_t}). \]
Evaluating the above expression at $n+1$ linearly independent values of
$\av$ allows us to conjecture $L$ and $M$.  Furthermore, the data 
revealed that inevitably $L<M$, so $q^L$ was factored out front of the expression,
and we renamed $q^{M-L}$ by $q^M$.
Conjecture~\ref{q11m1m1} was obtained similarly.

\section{Status of the conjectures}
As of this writing, the conjectures remain open.  Some twenty years ago,
J. Stembridge~\cite[p. 347, Corollary 7.4]{St}, in a different context, proved
that in the special case where $\av=\langle a,a,a,\dots, a \rangle$, and
$b_{\rho+1} = b_{\rho_2} = \cdots = b_{\rho_\tau} = -1$, for $\rho$ and $\tau$
satisfying $0\leqq \rho\leqq n$ and $1\leqq \tau\leqq n-\rho$,
\begin{equation}\label{stem}
 [ x_1^{b_1} x_2^{b_2} \cdots x_n^{b_n}] \mathcal F_n(\xv; \av; q) 
 = (-1)^\tau q^{b_1+b_2 + \cdots + b_\rho+ am} 
 \frac{ (q;q)_{an} (q^a;q^a)_\tau (q;q^a)_{\rho+\tau} }{ (q;q)_a^n (q;q^a)_n},
 \end{equation}
where $m=\rho\tau + \sum_{i=1}^\rho (i-1)b_i - \sum_{i=1}^{n-\rho-\tau} i b_{n-i+1}$.
One can check that Conjectures~\ref{q1m1}--\ref{q11m1m1} do in fact
agree with~\eqref{stem} in the instances where they overlap.

  It would of course be natural to investigate whether either or both proofs
of the $q$-Dyson conjecture~(\cite{ZB} and~\cite{GX}) could be adapted to
prove Conjectures~\ref{q1m1}--\ref{q11m1m1}.

\section{Possibilities for additional results}
It is likely that tractable formulas for additional co\"efficients in the Dyson
and $q$-Dyson products exist.  It seems quite plausible that the methods of this
paper would be sufficient for finding such formulas.   In particular, the next 
simplest case would likely be the co\"efficient of $x_r x_s / x_t^2$ in the
$q$-Dyson product, where $r$, $s$, and $t$ are distinct integers between
$1$ and $n$, with $n\geqq 3.$
\vskip 3mm
{\it Note added in proof:}  After the submission of this paper, Lv, Xin, and Zhou
announced a proof of Conjectures~\ref{q1m1}--\ref{q11m1m1} in ``A family of
$q$-Dyson style constant term identities," arXiv:0706.1009.

\bibliographystyle{amsalpha}

\begin{thebibliography}{SZ}

\bibitem[A]{A} G. E. Andrews, 
\textit{Problems and prospects for basic hypergeometric functions}, in: R. Askey (Ed.), 
\textit{The Theory and Application of Special Functions}, Academic Press, New York, 1975, 
pp. 191--224.

\bibitem[D1]{D1} F. J. Dyson, \textit{Statistical theory of the energy levels of complex systems I},
J. Math. Phys. \textbf{3} (1962) 140--156.

\bibitem[D2]{D2} F. J. Dyson, \textit{Disturbing the Universe}, Harper and Row, 1979.

\bibitem[Ga]{Ga} F. Garvan,  \textit{A $q$-product tutorial for a $q$-series MAPLE package}, The Andrews Festschrift (Maratea, 1998), Sm. Lothar. Combin. 42 (1999), Art. B42d, 27 pp. 

\bibitem[GX]{GX} I. M. Gessel and G. Xin, 
\textit{A short proof of the Zeilberger-Bressoud $q$-Dyson theorem}, Proc. Amer. Math. Soc. 
\textbf{134}
(2006) 2179--2187.

\bibitem[Go]{Go} I. J. Good, \textit{Short proof of a conjecture of Dyson}, J. Math. Phys. 
\textbf{11} (1970) 1884.

\bibitem[Gu]{Gu} J. Gunson, \textit{Proof of a conjecture of Dyson in the statistical theory of energy levels}, J. Math. Phys. \textbf{3} (1962) 752--753.

\bibitem[Si1]{Si1} A. V. Sills, \textit{Disturbing the Dyson Conjecture, in a
\emph{generally} GOOD way},
J. Combin.Theory Ser. A \textbf{113} (2006) 1368--1380.

\bibitem[Si2]{Si2} A. V. Sills, \texttt{qDysonConj} Maple package, available at\\
\texttt{http://www.georgiasouthern.edu/\~{}asills}.

\bibitem [SZ]{SZ} A. V. Sills and D. Zeilberger, \textit{Disturbing the Dyson Conjecture 
(in a GOOD Way)},
Experimental Math. \textbf{15} (2006) 187--191.

\bibitem[SZ2]{SZ2} A. V. Sills and D. Zeilberger, \texttt{GoodDyson} Maple package,
available at\\
\texttt{http://www.georgiasouthern.edu/\~{}asills} and \texttt{http://www.math.rutgers.edu/\~{}zeilberg}.

\bibitem[St]{St} J.R. Stembridge, 
\textit{First layer formulas for characters of $SL(n,\mathbb{C})$}, 
Trans. Amer. Math. Soc. \textbf{299} (1987) 319--350

\bibitem[W]{W} K. Wilson, 
\textit{Proof of a conjecture by Dyson}, J. Math. Phys. \textbf{3} (1962) 1040--1043.

\bibitem[ZB]{ZB} D. Zeilberger and D. M. Bressoud, 
\textit{A proof of Andrews' $q$-Dyson conjecture}, Discrete Math. \textbf{54} (1985) 201--224.

\end{thebibliography}

\end{document}